# Errata and supplements to: "Orthonormal RBF Wavelet and Ridgelet-like Series and Transforms for High-Dimensional Problems"

## W. Chen

*Department of Informatics, University of Oslo, P.O.Box 1080, Blindern, 0316 Oslo, Norway*

**Key words**: radial basis function, high-dimensional problems, fundamental solution, nonsingular general solution, wavelets, orthonormal RBF transformation and series.

## 1. Introduction

In recent years some attempts have been done to relate the RBF with wavelets [1,2] in handling high dimensional multiscale problems. To the author's knowledge, however, the orthonormal and bi-orthogonal RBF wavelets are still missing in the literature. By using the nonsingular general solution and singular fundamental solution of differential operator [3], recently the present author made some substantial headway to derive the orthonormal RBF wavelets series and transforms. The methodology can be generalized to create the RBF wavelets by means of the orthogonal convolution kernel function of various integral operators. In particular, it is stressed that the presented RBF wavelets does not apply the tensor product to handle multivariate problems at all.

This note is to correct some errata in ref. [3] and also to supply a few latest advances in the study of orthornormal RBF wavelet transforms.

## 2. Errata

Here we correct a few errors in the section 2 of ref. [3] involving orthonormal harmonic Bessel RBF series. To better understand, the readers are advised to refer to [3].

The formulas (4a,b) in ref. [3] is corrected as

$$\alpha_0 = \frac{2}{(n+2)S_n R^{n+1}} \int_{\Omega_R} r_{x\zeta}^{n/2} f(\zeta) d\Omega_\zeta \,,$$

$$\alpha_{jk} = \frac{2}{S_n R^{n+1} J_{n/2}(\lambda_j)^2} \left(\frac{\lambda_j}{2\pi}\right)^{1-n/2} \int_{\Omega_\zeta} r_{k\zeta}^{n/2} f(\zeta) J_{(n/2)-1}\left(\frac{\lambda_j r_{k\zeta}}{R}\right) d\Omega_\zeta, \quad j=1,2,\ldots, \quad k=1,2,\ldots,$$

where $S_n$ is the surface size of an unit $n$-dimensional sphere. Accordingly, the formula (5) in [3] should correctly be rewritten as

$$f(x) = \frac{1}{S_n R^{n+1}} \int_{\Omega_R} r_{x\zeta}^{n/2} f(\zeta) d\Omega_\zeta + \sum_{j=1}^{\infty} \sum_{k=1}^{\infty} \frac{2}{S_n R^{n+1} J_{n/2}(\lambda_j)^2} \int_{\Omega_R} r_{k\zeta}^{n/2} f(\zeta) J_{(n/2)-1}\left(\frac{\lambda_j r_{k\zeta}}{R}\right) d\Omega_\zeta r_{xk}^{1-n/2} J_{(n/2)-1}\left(\frac{\lambda_j r_{xk}}{R}\right)$$

## 3. Some updated results

In terms of ref. [3], some improved results are given below.

### 3.1. Bessel RBF transform

It is known [4] that

$$\lim_{j\to\infty}[\lambda_{n(j+1)} - \lambda_{nj}] = \pi \,,$$

So
$$\lim_{R\to\infty}\lim_{j\to\infty}\Delta\lambda_{nj}/R = \pi d\lambda.$$

Bessel function $J$ looks qualitatively like cosine waves when $x>>v$. The asymptotic form is

$$J_{(v)}(x) \approx \sqrt{\frac{2}{\pi x}}\cos(x-\frac{1}{2}v\pi-\frac{1}{4}\pi).$$

Assuming $f(x)$ absolutely integrable in $\Omega_\infty$, and approaching limit of $R$ and $\lambda$ in Eq. (5), we have the improved continuous Bessel transform

$$F(\lambda,\xi) = \int_{\Omega_\infty} r_{\xi\zeta}^{n/2} f(\zeta) J_{(n/2)-1}(\lambda r_{\xi\zeta}) d\Omega_\zeta,$$

$$f(x) = \frac{1}{C_{\varphi_n} S_n}\int_0^{+\infty}\int_{\Omega_\infty} F(\lambda,\xi) r_{x\xi}^{1-n/2} J_{(n/2)-1}(\lambda r_{x\xi})\lambda d\Omega_\xi d\lambda$$

### 3.2 Bi-orthogonal RBF transforms

The inverse bi-orthogonal RBF transform (16) in ref. [3] can be more precisely stated as

$$f(x) = C_g^{-1}\int_{-\infty}^{+\infty}\int_{\Omega_\infty} F(\lambda,\xi) g_n(\lambda r_{x\xi})\lambda d\Omega_\xi d\lambda,$$

where $C_g$ is the coefficient.

### 3.3. Harmonic time-space RBF transforms

By analogy with the Euclidean definition of distance variable, the time-space distance function was defined by Chen [5] for transient wave problem, namely,

$$\hat{r}_k = \sqrt{c^2\Delta t_k^2 - r_k^2},$$

where $r_k$ denotes the normal spatial Euclidean distance function, and $\Delta t_k = t - t_k$. $c$ is wave velocity. In dealing with time-dependent function, we have the time-space RBF wavelet approximation

$$f(x) = \alpha_0 + \sum_{j=1}^{\infty}\sum_{k=1}^{\infty}\alpha_{jk}\varphi_{nj}(\hat{r}_k) H(c\Delta t_k - r_k),$$

where $H$ is the Heaviside function. To discrete harmonic analysis, we still choose the non-singular general solutions of $n$-dimension Helmholtz equation as wavelet basis function. Similar to the procedure given in [3], we can determine the expansion coefficients by

$$\alpha_0 = \frac{2}{(n+2)S_n R^{n+1}}\int_{\Omega_R}\hat{r}_{x\zeta}^{n/2} f(\zeta) H(c\Delta t_{x\zeta} - r_{x\zeta}) d\Omega_\zeta,$$

$$\alpha_{jk} = \frac{2}{S_n R^{n+1} J_{n/2}(\lambda_j)^2}\left(\frac{\lambda_j}{2\pi}\right)^{1-n/2}\int_{\Omega_\zeta} r_{k\zeta}^{n/2} f(\zeta)$$
$$J_{(n/2)-1}\left(\frac{\lambda_j r_{k\zeta}}{R}\right) H(c\Delta t_{k\zeta} - r_{k\zeta}) d\Omega_\zeta,$$
$$j=1,2,\ldots,\quad k=1,2,\ldots,$$

The corresponding continuous Bessel and bi-orthogonal RBF transforms for transient wave function can be built by following the similar procedure given in [3].